\documentclass{article}
\usepackage{amsmath}
\usepackage{amsfonts}
\usepackage{amssymb}

\usepackage{graphicx}
\newcommand{\doble}{\mathop{\rightleftharpoons}}

\newcommand{\Z}{\mathbf Z}

\newcommand{\D}{\displaystyle}

\newtheorem{lemma}{\bf Lemma}{}{\rmfamily}
{}{\rmfamily}
{\itshape}{\rmfamily}
{}{\rmfamily}
\newtheorem{proposition}{\bf Proposition}{}{\rmfamily}
\newtheorem{algorithm}{\bf Algorithm}{}{\rmfamily}
\newtheorem{remark}{\bf Remark}{}{\rmfamily}
{}{\rmfamily}

\begin{document}

\title{On higher dimensional cocyclic Hadamard matrices} 
\author{V. \'{A}lvarez \and J.A. Armario \and M.D. Frau \and P. Real \thanks{All authors are partially supported by FEDER funds via the research projects FQM-296 and FQM-016 from JJAA.}} 







\maketitle

\begin{abstract}
Provided that a cohomological model for $G$ is known, we describe a method for constructing a basis for $n$-cocycles over
$G$, from which the whole set of $n$-dimensional $n$-cocyclic matrices over $G$ may be straightforwardly calculated. Focusing in the case
$n=2$ (which is of special interest, e.g. for looking for cocyclic Hadamard matrices), this method provides a basis for 2-cocycles in such a
way that representative $2$-cocycles are calculated all at once, so that there is no need to distinguish between inflation and
transgression 2-cocycles (as it has traditionally been the case until now). When $n>2$, this method provides an uniform way of looking for higher dimensional {\em $n$-cocyclic} Hadamard matrices for the first time. We illustrate the method with some examples, for $n=2,3$. In particular, we give some examples of improper 3-dimensional $3$-cocyclic Hadamard matrices.


\end{abstract}

\section{Introduction}

Hadamard matrices are square matrices with entries $\pm 1$ such that their rows are pairwise orthogonal. They were first noticed by Sylvester in problems related to tessellated pavements, and later on by Hadamard related to the maximal determinant problem, which asks for the largest determinant for all matrices  with entries $\pm 1$. Recommended references on Hadamard matrices and
their applications are \cite{HW78} and more recently \cite{Hor07}.

It is easy to prove that the size
of Hadamard matrices must be 1, 2 or a multiple of 4. Nevertheless, it is an open question whether Hadamard matrices exist for every size $4t$. This is known as the {\em Hadamard Conjecture}. A lot of work has been made concerning this conjecture, and the ways in which Hadamard matrices might be constructed. All known  construction methods such as Sylvester Hadamard matrices, Paley Hadamard matrices, Williamson Hadamard matrices or Ito Hadamard matrices fail to yield Hadamard matrices for every order which is a multiple of 4. The most promising methods which may lead to a constructive proof of the Hadamard Conjecture are the two-circulant core construction \cite{Kot13} and the cocyclic approach \cite{HdL95}. 

A matrix $M$ is said to be {\em cocyclic} if there exist a group $G$ of order $|G|=4t$ and a 2-cocycle
$\psi :G \times G \rightarrow \{ \pm 1\}$ such that $M=(\psi (g_i,g_j))$. The main advantages of the cocyclic framework concerning Hadamard matrices  may be summarized in the following facts:

\begin{itemize}

\item Determining if a $\pm 1$ matrix
$M$ of order $4t$ is Hadamard would consist of checking whether the rows of $M$ are pairwise orthogonal, which would require in total $O(t^3)$ operations. In fact, since these matrices give the solution to the problem of the maximal determinant problem of matrices of size $4t \times 4t$ \cite{Had93}, it suffices to check whether $|H|=(4t)^{2t}$, which actually can be done in less than $O(t^3)$ time, using some fast matrix multiplication techniques.

Nevertheless, a cocyclic matrix $M=(\psi (g_i,g_j))$ is
Hadamard if and only if the summation of each row but the first is zero (see the cocyclic Hadamard test of \cite{HdL95}), which requires at most $O(t^2)$ operations.

\item The search space is reduced to the set of cocyclic matrices over a given group, instead of the whole set of matrices with entries in $\{-1,1\}$. In spite of this fact, the proportion of cocyclic Hadamard matrices among the full set of cocyclic matrices seems to be not significantly different from that of usual Hadamard matrices among $\pm 1$ matrices, as calculations in \cite{AAFG12} suggest.   \end{itemize}

In this paper we are concerned with higher dimensional Hadamard matrices, and the way in which the cocyclic framework may be used in this context. For commodity, in what follows any 3D-matrix $A(a_{i_1,i_2,i_3})_{1\leq i_j \leq v,1 \leq j \leq 3}$ will be described by listing its 2-dimensional horizontal sections $A(a_{i_1,i_2,i})$, for $1 \leq i \leq v$. As usual, negative entries will be denoted simply by $-$.

As introduced by Shlichta in \cite{Shl71,Shl79}, an {\em improper} $n$-dimensional Hadamard matrix of order $v$ is a $(\pm 1)$ array $A=(a_{i_1,i_2,\ldots,i_n})_{1\leq i_j \leq v,1 \leq j \leq n}$ such that all its parallel $(n-1)$-dimensional sections are mutually orthogonal; that is,
for each $1 \leq l \leq n$, and for all indices $x$ and $y$ in dimension $l$, $$ \sum _{j  \neq l} \sum _{1 \leq i_j \leq v} a_{i_1, \ldots , x, \ldots, i_j, \ldots, i_n} \cdot a_{i_i, \ldots , y, \ldots , i_j, \ldots , i_n}=v^{n-1} \delta _{xy}.$$

For instance, the following improper 3-dimensional Hadamard matrix  $A=(a_{i_1,i_2,i_3})_{1 \leq i_j\leq 2}$ may be found in \cite{Shl79}, where $A(i_1,i_2,1)=\left( \begin{array}{cc} 1&1\\-&1 \end{array} \right)$ and $A(i_1,i_2,2)=\left( \begin{array}{cc} -&-\\-&1 \end{array} \right)$. Notice that the section $A(a_{1,i_2,i_3})= \left( \begin{array}{cc} 1&1\\ 0&0 \end{array} \right)$ is not Hadamard, in the usual sense.

Notice that an improper $n$-dimensional Hadamard matrix may have stronger orthogonality properties in some dimensions. For instance, an $n$-dimensional Hadamard matrix $A$ is  termed {\em proper} (see \cite{Hor07,Yan01} for details) if any parallel rows in $A$ are orthogonal; that is, for each pair of dimensions $j,l$, for all indices $x$ and $y$ in dimension $l$, and for each set of fixed indices in the other $n-2$ dimensions, $$  \sum _{1 \leq i_j \leq v} a_{i_1, \ldots , x, \ldots, i_j, \ldots, i_n} \cdot a_{i_i, \ldots , y, \ldots , i_j, \ldots , i_n}=v \delta _{xy}.$$

The matrix of the example above is an improper not proper 3-dimensional Hadamard matrix. Nevertheless, in \cite{Shl79} one may find some proper 3-dimensional Hadamard matrices, such as $A=(a_{i_1,i_2,i_3})_{1 \leq i_j\leq 2}$, consisting of the horizontal sections $A(i_1,i_2,1)=\left( \begin{array}{cc} 1&-\\1&1 \end{array} \right)$ and $A(i_1,i_2,2)=\left( \begin{array}{cc} 1&1\\-&1 \end{array} \right)$.

There are some well-known methods for constructing both improper and proper higher-dimensional Hadamard matrices. The interested reader is referred to \cite{Hor07,Yan01,dLF11} and the references there cited.

It follows that every planar section of a proper $n$-dimensional Hadamard matrix of order $v$ is a Hadamard matrix of order $v$ itself, so $v$ must be 2 or a multiple of 4. Actually, it is known that proper  $n$-dimensional Hadamard matrices of order $v=4t$ exist if and only if usual Hadamard matrices of order $v=4t$ do exist (see \cite{dLa86,Yan86} for details).


It is worthwhile that although planar (and therefore proper $n$-dimensional) Hadamard matrices can only exist on orders $v >2$ multiple of 4, improper higher-dimensional Hadamard matrices may exist on even orders not multiple of 4 (see \cite{Yan01}). It is easy to prove, though, that the order of an improper higher-dimensional Hadamard matrix must be even, in any case.

Few is known about the existence of improper $n$-dimensional Hadamard matrices of order $v$. There are well known examples of improper $3,4$-dimensional Hadamard matrices of order 6 (see \cite{Yan01}). Anyway, it is not known whether improper $n$-dimensional Hadamard matrices exist for all even orders $v=2t$, $n \geq 3$, for odd $t$. Furthermore, only $3$-dimensional improper Hadamard matrices of order $2 \cdot 3^b$ are known, $b\geq 1$. The interested reader is referred to \cite{Yan01} for more details.

However, despite its potential, the subject of higher-dimensional Hadamard matrices remains seriously under-developed. The aim of this paper is to establish a basis for the study of higher-dimensional Hadamard matrices from the cocyclic point of view.

As introduced in \cite{Shl71,Shl79}, the notion of higher dimensional Hadamard matrix is a natural generalization of usual Hadamard matrices to the case of $n$-dimensional matrices, for $n\geq 2$. One could wonder whether the cocyclic framework could also be taken into account for higher dimensional Hadamard matrices.

A precedent in the literature on this subject may be found in de Launey's work in \cite{dLa89} (later extended in \cite{dLH93,dLF11}). Here,  a way to construct a proper $n$-dimensional Hadamard matrix $A=(a_{i_1,\ldots,i_n})$ from a 2-cocyclic planar Hadamard matrix $H=(\psi (g_i,g_j))$ is described, so that  \begin{equation}\label{delauney} a_{i_1,\ldots,i_n}= \prod _{k=2}^n \psi (\prod _{j=1}^{k-1}g_{i_j},g_{i_k}), \end{equation} for $\psi:G \times G \rightarrow \{-1,1\}$ being a 2-cocycle over $G$.

We wonder if one could go farther and use $n$-cocycles $\phi: G \times \stackrel{n}{\ldots} \times G \rightarrow \{ 1,-1\}$ in an attempt to find $n$-dimensional Hadamard matrices $A=(\phi(i_1,\ldots,i_n))$, for $n>2$. If so, this would extend both the works of Shlichta and de Launey described before. Notice that given a $2$-cocycle $\psi: G \times G \rightarrow \{1,-1\}$, a function $f_{\psi}: G \times \stackrel{n}{\ldots}\times G \rightarrow \{1,-1\}$ defined as in (\ref{delauney}), so that $\D f_{\psi}(i_1,\ldots,i_n)= \prod _{k=2}^n \psi (\prod _{j=1}^{k-1}g_{i_j},g_{i_k})$, does not define a $n$-cocycle $\phi$ in general.

Recall that for a map $\phi: G\times \stackrel{n}{\ldots}\times G \rightarrow \{1,-1\}$ being a $n$-cocycle over $G$, it is necessary and sufficient that for all $g_{i_1},\ldots,g_{i_{n+1}}\in G,$ \begin{equation}\label{3cocyclic} \phi(g_{i_2},\ldots,g_{i_{n+1}})\cdot \phi(g_{i_1},\ldots,g_{i_{n}}) \cdot\prod_{k=1}^{n}\phi (g_{i_1},\ldots,g_{i_k}g_{i_{k+1}},\ldots,g_{i_{n+1}})=1. \end{equation}
In particular, consider $G=\Z _2$ and the 2-cocycle $\psi: \Z_2 \times \Z_2 \rightarrow \{1,-1\}$, so that $\psi(0,0)=\psi(0,1)=\psi(1,0)=1$, $\psi (1,1)=-1$. By $(1)$, the matrix $A=(f_{\psi}(i_1,i_2,i_3))$ consisting of $f_{\psi}(i_1,i_2,i_3)=\psi(i_1,i_2)\cdot \psi(i_1i_2,i_3)$ defines a proper 3-dimensional Hadamard matrix, whose horizontal sections are $A(i_1,i_2,0)=\left( \begin{array}{cc} 1&1\\1&-\end{array}\right)$ and $A(i_1,i_2,1)=\left( \begin{array}{cc} 1&-\\-&-\end{array}\right)$. However, $f_{\psi}$ defines by no means  a $3$-cocycle over $\Z _2$, since substituting  $(g_{i_1},g_{i_2},g_{i_3},g_{i_4})=(0,0,1,1)$ in (\ref{3cocyclic}) leads to $$f_{\psi}(0,1,1)\cdot f_{\psi}(0,0,1)\cdot f_{\psi}(0,1,1)\cdot f_{\psi}(0,1,1)\cdot f_{\psi}(0,0,0)=-1 \neq 1.$$ Thus it makes sense going beyond de Launey's work about constructing $n$-dimensional proper Hadamard matrices from orthogonal $2$-cocycles, and trying to look for $n$-dimensional Hadamard matrices coming from $n$-cocycles, for $n>2$.

In these circumstances, two main questions should be studied, for $n>2$:

\begin{enumerate}

\item Is there any cocyclic test for $n$-cocyclic $n$-dimensional (proper/improper) Hadamard matrices?

\item Is there any effective way to construct $n$-cocyclic $n$-dimensional matrices?

\end{enumerate}

Although we have studied the first of these problems, unfortunately we have not been able to isolate any simpler characterization for a $n$-cocyclic $n$-dimensional matrix to be (proper/improper) Hadamard for the moment, for $n>2$.

In order to answer the second question, though, a prerequisite   is to determine a basis for $n$-cocyclic matrices over a given finite group $G$. Here we  provide a method for constructing $n$-cocyclic $n$-dimensional matrices over some finite groups $G$. Nevertheless, the difficult step is not just constructing such matrices, but  determining which among them satisfy the proper/improper Hadamard condition.

Focusing in the case $n=2$ (which is of special interest, e.g. for looking for cocyclic Hadamard matrices), the Universal Coefficient
Theorem provides a decomposition of representative 2-cocycles over a group $G$ as the direct sum of the inflation and transgression
cocycles over $G$, $$H^2(G,\Z_2)\cong Ext(G/[G,G],\Z_2)\oplus Hom(H_2(G),\Z_2).$$

Until now, all the methods which look for a basis for 2-cocycles uses this decomposition, so that two different processes have to be
performed.

In spite of this fact, there is a chance to calculate a basis for $H^2(G,\Z_2)$ all at once in a straightforward manner, provided that a
cohomological model for $G$ is known, so that there is no need to distinguish between inflation and transgression 2-cocycles, as it has
traditionally been the case. It is the cohomological analog of the homological reduction method described in \cite{AAFR08a}. An early
precedent of this technique is located in \cite{GL00}, focusing on $p$-groups.

This procedure may be extended in order to construct a basis for $n$-cocycles as well as $n$-dimensional $n$-cocyclic matrices over a group
$G$ for which a cohomological model $cohG$ is known. We call this process the ``cohomological reduction method''.

The term {\em cohomological model} refers to a special type of homotopy equivalence (termed {\em contraction} \cite{EM54}) $$\D^{\phi:}{
Hom(\bar{B}(\Z[G]),\Z_2)} \doble_K^F cohG$$ from the set of homomorphisms of the {\em reduced bar construction} (i.e. the reduced complex
associated to the standard bar resolution \cite{EM54}) of the group $G$ onto $\Z$, to a differential graded comodule of finite type
$cohG$. Thus $$H^*(G)=H^*(Hom(\bar{B}(\Z[G]),\Z_2))=H^*(cohG)$$ and the $n$-cohomology of $G$ and its representative $n$-cocycles may be
effectively computed from those of $cohG$, by means of Veblen's algorithm \cite{Veb31} (involving the Smith's normal forms of the matrices
representing the codifferential operator).

In particular, if $\D^{\phi:}{ \bar{B}(\Z[G])} \doble_K^F (hG,d)$ defines a homological model for $G$, then $$\D^{\phi^*:}{
Hom(\bar{B}(\Z[G]),\Z_2)} \doble_{F^*}^{K^*} (Hom(hG,\Z_2),d^*)$$ defines a cohomological model for $G$. This way, the set of groups for
which a cohomological model is known includes those groups for which a homological model is known. Consequently, this method extends the
homological reduction method described in \cite{AAFR08a}.

The paper is organized as follows.  In section 2 we describe the cohomological reduction method itself, that is, how to construct a full basis for
$n$-cocycles over $G$ from a cohomological model $cohG$ for $G$. Section 3 is devoted to show several 2 and 3 dimensional examples, including the
well-known cases of dihedral groups $D_{2t}$ and abelian groups $\Z_{2t}$, $\Z_{2t}\times \Z_2$ and $\Z_t \times \Z_2^2$ for clarity. All the calculations have
been made with aid of some packages in Mathematica provided by the authors in \cite{AAFR06h,AAFR06i,AAFR06j}. Some conclusions and future work is described in the last section.

\section{Describing the cohomological reduction method}






The cohomological reduction method provides a computationally efficient way to lift  the cohomological information from a cohomological model
$$\D^{\phi:}{ Hom(\bar{B}(\Z[G]),\Z_2)} \doble_K^F cohG$$ for a group $G$ to the group itself. The injection morphism  $$K:cohG \longrightarrow
Hom(\bar{B}(\Z[G]),\Z_2)$$ helps in this task.

Let ${\cal B}_{i-1} =\{{\bf u_1},\ldots,{\bf u_q}\}$, ${\cal B}_{i} =\{{\bf e_1},\ldots,{\bf e_r}\}$ and ${\cal B}_{i+1} =\{{\bf v_1},\ldots,{\bf
v_s}\}$ be the corresponding basis for $cohG$ at dimensions $i-1$, $i$ and $i+1$, respectively. Attending to Veblen's algorithm, since
$H^i(G;\Z_2)\cong H^i(cohG;\Z_2)=Ker\; d^i/ Im \; d^{i-1},$ we need to calculate the binary (i.e. with coefficients in $\Z_2$) Smith Normal Forms of
the matrices representing the codifferential operators $d^{i-1}$ and $d^i$, $$M_{i-1}(d)=\left(
\begin{array}{c}
  d({\bf u_1}) \\
      \vdots\\
  d({\bf u_q})
\end{array}\right)_{q\times r}\hspace*{-.25cm}D_{i-1}=\left( \begin{array}{c|c}
  I_l &0 \\
    \hline 0& 0 \end{array}\right)_{q\times r} \hspace*{-.25cm}M_{i}(d)=\left(
\begin{array}{c}
  d({\bf e_1}) \\
      \vdots\\
  d({\bf e_r})
\end{array}\right)_{r\times s}\hspace*{-.25cm} D_i=\left( \begin{array}{c|c}
  I_k &0 \\
    \hline 0& 0 \end{array}\right)_{r\times s}$$ so that $H^i(G;\Z_2)\cong
H^i(cohG;\Z_2)\cong \Z _2^{r-k-l}$.

Furthermore, some change of basis matrices $P_j$ and $Q_j$ exist, for $j=i-1,i$, such that 
$$\begin{array}{ccc} {\cal B}_j &\stackrel{M_j(d)}{\longrightarrow} &{\cal B}_{j+1} \\
\makebox[0pt][r]{$\scriptstyle P_j$}{\uparrow}  & \mbox{ }^\# &\downarrow
\makebox[0pt][l]{$\scriptstyle Q_j$}\\
\bar{\cal B}_j &\stackrel{D_j}{\longrightarrow} &\bar{\cal B}_{j+1} \\
\end{array} \qquad \qquad D_j=P_j\cdot M_j(d)\cdot Q_j $$

Now we proceed according to the following steps:
\begin{enumerate}

\item A basis $\bar{\cal C}$ for $Im\; d^{i-1}$ with regards to ${\cal B}_i$ is obtained from the first $l$ columns of $Q_{i-1}^{-1}$. Thus
$\bar{\cal C}$ is a basis for $i$-coboundaries over $cohG$.

\item A basis $\bar{\cal D}$ for $Ker\; d^i$ with regards to ${\cal B}_i$ is obtained from the $r-k$ last rows of $P_i$. Thus $\bar{\cal D}$ is a
basis for $i$-cocycles over $cohG$.

\item Select those $r-k-l$ elements in $\bar{\cal D}$ which are not linear combinations of the elements in $\bar{\cal C}$. This forms a basis
$\bar{\cal B}$ for representative $i$-cocycles over $cohG$.

\item Lift $\bar{\cal B}$ to the correspondent basis ${\cal B}$ in $H^i(\bar{B}_i(\Z[G]);\Z_2)$ by means of the injection $K$.

\end{enumerate}

Graphically, 
$$\begin{array}{ccc}Hom(\bar{B}_i(\Z[G]),\Z_2)&\stackrel{K}{\longleftarrow}& {\cal B}_i \\
  &  &\uparrow
\makebox[0pt][l]{$\scriptstyle P_i,Q_{i-1}^{-1}$}\\ & &\bar{\cal B}_i \\
\end{array}$$

\begin{proposition}\label{basecohom} The scheme above defines a basis ${\cal B}$ for representative $i$-cocycles over $G$. \end{proposition}


A basis for $i$-coboundaries may be obtained by Linear Algebra. More concretely, denote $\partial_{[g_1,\ldots,g_{i-1}]}:\bar{B}_i(\Z[G])
\rightarrow \Z_2$ the $i$-coboundary associated to the characteristic map $\delta_{[g_1,\ldots,g_{i-1}]}$ of the element $[g_1,\ldots,g_{i-1}]\in
\bar{B}_{i-1}(\Z[G])$, $$\partial_{[g_1,\ldots,g_{i-1}]}([h_1,\ldots,h_i])=
\delta_{[g_1,\ldots,g_{i-1}]}([h_2,\ldots,h_i])+\delta_{[g_1,\ldots,g_{i-1}]}([h_1,\ldots,h_{i-1}])+$$ $$ \sum
_{j=1}^{i-1}\delta_{[g_1,\ldots,g_{i-1}]}([h_1,\ldots,h_jh_{j+1},\ldots,h_i])\qquad mod \; 2$$ Take the $4t\times \stackrel{i}{\ldots} \times 4t$
$i$-dimensional matrix $M_{\partial_{[g_1,\ldots,g_{i-1}]}}$ related to $\partial_{[g_1,\ldots,g_{i-1}]}$ as vectors of length $4^it^i$. Moreover,
consider the $4^{i-1}t^{i-1}\times 4^it^i$ matrix $C$ whose rows are the vectors $M_{\partial_{[g_1,\ldots,g_{i-1}]}}$. Then a row reduction on $C$
leads to a basis for i-coboundaries. It suffices to keep track of those coboundaries $\partial_{[g_1,\ldots,g_{i-1}]}$ whose transformed rows in
$M_{\partial_{[g_1,\ldots,g_{i-1}]}}$ after the row reduction are not zero.

\begin{lemma} \label{basecobordes} The morphisms $\partial_{[g_1,\ldots,g_{i-1}]}$ above define a basis for i-coboundaries. \end{lemma}




The cohomological reduction method provides then the following algorithm for computing $n$-dimensional $n$-cocyclic matrices over $G$.

\begin{algorithm} \label{mrh} {\bf (cohomological reduction method)}

{\tt \noindent {\sc Input:} group with cohomological model $\{G,cohG,F,K,\phi\}$ \\
\\
\hspace*{.25cm}  Construct a basis ${\cal E}$ for $i$-coboundaries (Lemma \ref{basecobordes}).\\
\hspace*{.25cm}  Construct a basis ${\cal B}$ for representative $i$-cocycles (Proposition \ref{basecohom}).\\
\\
\noindent {\sc Output:} By juxtaposition, a basis  ${\cal B}\cup {\cal E}$ for $i$-cocycles over $G$.} \end{algorithm}


\section{Examples}

We next show how the cohomological reduction method works for constructing basis for 2-cocycles and 3-cocycles over some groups (which have been
shown to provide many 2-cocyclic Hadamard matrices, see \cite{AAFR07}).

All the executions and examples of this section have been worked out with aid of the {\em Mathematica 4.0} notebooks \cite{AAFR06h,AAFR06i} described
in \cite{AAFR06c,AAFR06g} (for constructing homological models) and \cite{AAFR06b} (in order to form a basis for 2-cocycles from which the search for
cocyclic Hadamard matrices is then developed), running on a {\em Intel(R) Core(TM) i3 CPU, M330, 2.13GHz, 4,00GB RAM, 64 bits}.

In the sequel, the elements of a product $A\times B$ are ordered as the rows of a matrix indexed in $|A|\times |B|$. For instance, if $|A|=r$ and
$|B|=c$, the ordering is $$\langle a_1b_1,a_1b_2,\ldots, a_1b_c,a_2b_1,a_2b_2,\ldots,a_2b_c,\ldots,a_rb_1,\ldots, a_rb_c \rangle$$ The elements in
the group $G$ are labeled from 1 to $|G|$, accordingly to this ordering.

The back negacyclic matrix of order $j$ is denoted by $BN_j=\left( \begin{array}{crcr} 1&1&\cdots & 1\\
1&&_{\cdot}\cdot^{\cdot}&-1\\ \vdots&_{\cdot}\cdot^{\cdot}&_{\cdot}\cdot^{\cdot}&\vdots\\ 1&-1&\cdots&-1\end{array} \right)_{j\times j}$, as usual. The square matrix of order
$n$ formed all of 1s is denoted by $1_n$. The Kronecker product of matrices is denoted by $\otimes$, so that $A \otimes B$ is the block matrix
$\left( \begin{array}{ccc} a_{11}B & \ldots & a_{1n}B\\ \vdots && \vdots\\ a_{n1}B & \ldots & a_{nn}B \end{array}\right)$. The Hadamard (pointwise)
product of matrices is simply denoted as $A\cdot B$.  We use the Kronecker-Iverson notation $[b]$ (see \cite{GKP89}), which evaluates to 1 for
Boolean expressions $b$ having value true, and to 0 for those having value false. Finally, the notation $[x]_m$ refers to $x \; mod \; m$.

Let consider the families of groups below (assume $\Z_k=\{0,1,\ldots,k-1\}$ with additive law).

\begin{enumerate}

\item $G_1^t=\Z_{2t} \times \Z_2$.


\item $G_2^t=\Z_t \times \Z_{2}^2=\Z_t \times(\Z_2 \times \Z _2)$. Notice that $G_2^t\simeq G_1^t$ for odd $t$.


\item $G_3^t=\Z_{2t}$.

\item $G_4^t=D_{2t}$, for odd $t$.

\end{enumerate}

In this section we will construct a cohomological model for $G_i^t$ from the homological models for $G_i^t$ described in \cite{AAFR07}, so that if
$\D^{\phi:}{ \bar{B}(\Z[G_i^t])} \doble_K^F (hG_i^t,d)$ defines a homological model for $G_i^t$, then $$\D^{\phi^*:}{ Hom(\bar{B}(\Z[G_i^t]),\Z_2)}
\doble_{F^*}^{K^*} (Hom(hG_i^t,\Z_2),d^*)$$ defines a cohomological model for $G_i^t$. Here, as usual, we use $-^*$ for noting the dual object for
$-$.

More concretely, it suffices to take duals on the basis ${\cal B}_k$ for $hG_i^t$ on degree $1\leq k \leq 4$, and the differential operators
$d_j:B_{j+1}\rightarrow B_j$ and the projections  $F_j:\bar{B}_j(\Z[G_i^t]) \rightarrow B_j$ for $2\leq j \leq 3$. Recall that
$$\bar{B}_k(\Z[G])=\langle [g_1,\ldots,g_k]:\; g_j\in G \rangle.$$

Notice that the matrices $P$ and $Q$ involved in the calculation of the Smith Normal Form, $D$, for a matrix $A$ (so that $D=P\cdot A\cdot Q$) are
not uniquely determined, in general. In the sequel we will use the matrices coming from the SmithNormalForm package programmed in \cite{AAFR06j}.

For illustrating the method, we will display explicitly the computation for the case of $G_1^t=\Z_{2t} \times \Z_2$, giving a basis for 2- and 3- cocycles over $G_1^t$ in detail. On the contrary, just a summary of the results of the analog computations for $G_2^t$ will be presented. We will show that the basis of 2-cocycles obtained so far and the basis for 2-cocycles calculated in \cite{AAFR07} are equivalent. Later on, we will use the basis for 3-cocycles in order to look for $3$-dimensional Hadamard matrices of order $4t$ over either $G_1^t$ or $G_2^t$.

In addition, we will calculate some basis for 3-cocycles over $G_3^t$ and $G_4^t$, from which a search for 3-dimensional improper Hadamard matrices of order $2t$ will be performed. Surprisingly, although $D_{4t}$ is prolific giving rise to many 2-cocyclic Hadamard matrices, $G_4^t=D_{2t}$ will show to be not suitable for looking for 3-dimensional 3-cocyclic improper Hadamard matrices, since there are not any representative 3-cocycles over $G_4^t$, and hence the set of 3-cocycles over $D_{2t}$ reduces to the set of 3-coboundaries over $D_{2t}$.

\subsection{Basis for 2-,3-cocycles over $G_1^t=\Z_{2t} \times \Z_2$}$\mbox{ }$\\

Notice that the $i$-th element of $G_1^t$ corresponds to $(\lfloor \frac{i-1}{2}\rfloor,[i-1]_2)\in \Z_{2t} \times \Z_2$. Conversely, the element
$(i_1,i_2)\in \Z_{2t} \times \Z_2$ corresponds to the $i$-th element of $G_1^t$, for $i=2i_1+i_2+1$.

A homological model for $G_1^t$ is described in \cite{AAFR07} and consists of

\noindent ${\cal B}_1=\{u_1,u_2 \}$, ${\cal B}_2=\{e_1,e_2,e_3 \}$, ${\cal B}_3=\{v_1,v_2,v_3,v_4 \}$,\\ $d_2(e_1)=2t\cdot u_1$, $d_2(e_3)=2\cdot
u_2$, $d_3(v_2)=2t\cdot e_2$, $d_3(v_3)=-2\cdot e_2$,\\ $F_1[i]=i_1\cdot u_1+i_2\cdot u_2$, $F_2[i|j]=[i_1+j_1\geq 2t]\cdot e_1+ i_1j_2\cdot
e_2+[i_2+j_2\geq 2]\cdot e_3$.

In these circumstances, it may be checked that

\noindent ${\cal B}_1^*=\{u_1^*,u_2^* \}$, ${\cal B}_2^*=\{e_1^*,e_2^*,e_3^* \}$, ${\cal B}_3^*=\{v_1^*,v_2^*,v_3^*,v_4^* \}$,\\
$d^1=d_2^*=0$, $d^2=d_3^*=0$,\\ $F^*(e_1^*)([i|j])=e_1^*(F([i|j]))=\left[ [i_1+j_1\geq 2t]\right]_2$,
$F^*(e_2^*)([i|j])=e_2^*(F([i|j]))=\left[i_1j_2\right]_2$, $F^*(e_3^*)([i|j])=e_3^*(F([i|j]))=\left[[i_2+j_2\geq 2]\right]_2.$

From these data, it may be checked that $\begin{array}{|c|c|c|c|} \hline D_1&Q^{-1}_1&D_2&P_2\\ \hline  0_{2\times 3} & I_3 & 0_{3\times 4} & I_3 \\
\hline
\end{array}$

Thus $H^2(cohG_1^t;\Z_2)=\Z_2^3$ and ${\cal B}_2^*=\{e_1^*,e_2^*,e_3^* \}$ is a basis for representative 2-cocycles over $cohG_1^t$. Accordingly, a
basis for representative 2-cocycles over $G_1^t$ is given by $\{F^*(e_1^*),F^*(e_2^*),F^*(e_3^*)\}=\{BN_{2t}\otimes 1_2,1_t \otimes \left(
\begin{array}{cccc}1&1&1&1\\1&1&1&1\\1&-&1&-\\1&-&1&- \end{array}\right),1_{2t}\otimes
BN_2\}$.

This basis may be extended to a full basis for 2-cocycles over $G_1^t$, by simply juxtaposing a basis for 2-coboundaries over $G_1^t$ (see Lemma
\ref{basecobordes}). This basis is $\langle
\partial_2,\ldots,\partial_{4t-2}\rangle$ as it was pointed out in \cite{AAFR07}.

\begin{remark} A basis $\{\beta_1,\beta_2,\gamma_1\}$ for representative 2-cocycles over $G_1^t$ was already determined in \cite{AAFR07}
(notice that $G_1^t$ was denoted $G_2^t$ there). It may be checked that $\beta_1=F^*(e_3^*)$ and $\gamma_1=F^*(e_2^*)$. If $2t=2^rq$, for
odd $q$, then $$\D \beta_2=1_q \otimes BN_{2^r}\otimes 1_2=F^*(e_1^*)\cdot \prod _{k=0}^{\lfloor \frac{q}{2}\rfloor-1} \partial
_{k2^{r+2}+2^{r+1}+1} \ldots \partial _{k2^{r+2}+2^{r+2}}.$$ \end{remark}

Now we compute a basis for 3-cocycles over $G_1^t$. From calculations in \cite{AAFR07,AAFR06c}, it is easy to derive that the homological model described above may be extended to degree 4, so that ${\cal B}_4=\{w_1,w_2,w_3,w_4,w_5 \}$,
$d_4(w_1)=2t\cdot v_1$, $d_4(w_3)=2t\cdot v_3+2\cdot v_2$, $d_4(w_5)=2\cdot v_4$, $F[i|j|k]=[ k_1[i_1+j_1\geq 2t]\cdot v_1+ k_2[i_1+j_1\geq 2t]\cdot v_2+i_1[j_2+k_2\geq 2]\cdot v_3+ k_2[i_2+j_2\geq
2]\cdot v_4]_2$.

Consequently, we consider ${\cal B}_4^*=\{w_1^*,w_2^*,w_3^*,w_4^*,w_5^* \}$, $d^3=0$, $F^*(v_1^*)([i|j|k])=v_1^*(F([i|j|k]))=[ k_1[i_1+j_1\geq 2t]]_2$,
$F^*(v_2^*)([i|j|k])=v_2^*(F([i|j|k]))=[ k_2[i_1+j_1\geq 2t]]_2$, $F^*(v_3^*)([i|j|k])=v_3^*(F([i|j|k]))=[i_1[j_2+k_2\geq 2]]_2$, $F^*(v_4^*)([i|j|k])=v_4^*(F([i|j|k]))=[k_2[i_2+j_2\geq
2]]_2$.

From these data, it may be checked that $\begin{array}{|c|c|c|c|} \hline D_2&Q^{-1}_3&D_3&P_3\\ \hline  0_{3\times 4} & I_4 & 0_{4\times 5} & I_4 \\
\hline
\end{array}$

Thus $H^3(cohG_1^t;\Z_2)=\Z_2^4$ and ${\cal B}_3^*=\{v_1^*,v_2^*,v_3^*,v_4^* \}$ is a basis for representative 3-cocycles over $cohG_1^t$.
Accordingly, a basis for representative 3-cocycles over $G_1^t$ is given by $\{F^*(v_1^*),F^*(v_2^*),F^*(v_3^*),F^*(v_4^*)\}$.

It may be straightforwardly checked that the $F^*(v_i^*)$ are the 3D-matrices whose horizontal sections are given by:

 \begin{itemize}

 \item $F^*(v_1^*)$: $\{J_{4t},J_{4t},BN_{2t} \otimes 1_2,BN_{2t} \otimes 1_2, \stackrel{t}{\ldots},J_{4t},J_{4t},BN_{2t} \otimes 1_2,BN_{2t} \otimes 1_2\}$.

 \item $F^*(v_2^*)$: $\{J_{4t},BN_{2t} \otimes 1_2, \stackrel{2t}{\ldots},J_{4t},BN_{2t} \otimes 1_2\}$.

 \item $F^*(v_3^*)$: $\{J_{4t},1_{t}\otimes A, \stackrel{2t}{\ldots},J_{4t},1_{t} \otimes A\}$, for $A= \left( \begin{array}{cccc} 1&1&1&1\\ 1&1&1&1\\ 1&-&1&- \\ 1&-&1&- \end{array} \right)$.

  \item $F^*(v_4^*)$: $\{J_{4t},1_{2t}\otimes BN_{2}, \stackrel{2t}{\ldots},J_{4t},1_{2t}\otimes BN_{2}\}$.

\end{itemize}

This basis may be extended to a full basis for 3-cocycles over $G_1^t$, by simply juxtaposing a basis for 3-coboundaries over $G_1^t$ (see Lemma
\ref{basecobordes}). This basis is given by $\langle
\partial_1,\ldots,\partial_{16t^2-4t-2},\partial_{16t^2-4t+1}\rangle$ by direct inspection.

\subsection{Basis for 2-,3-cocycles over $G_2^t=\Z_t \times \Z_{2}^2=\Z_t \times(\Z_2 \times \Z _2)$}$\mbox{ }$\\

Progressing from the homological model for $G_2^t$ described in \cite{AAFR07}, the cohomological reduction method straightforwardly provides that $H^2(cohG_2^t;\Z_2)=\left\{ \begin{array}{ll} \Z_2^3 & \mbox{ for $t$ odd,}\\  \Z_2^6 & \mbox{ for $t$ even.} \end{array}\right.$ Depending on whether $t$ is odd or even, a basis for representative $2$-cocycles over $G_2^t$ is obtained considering just the first three elements or the full set of six matrices of the following set:
$\{BN_t \otimes 1_4,1_{\frac{t}{2}} \otimes K_2,1_{\frac{t}{2}} \otimes K_3,1_t \otimes BN_2 \otimes 1_2,1_t \otimes K_1,1_{2t}\otimes BN_2\}$ where
the matrices $K_1,K_2,K_3$ are given by $$\begin{array}{|c|c|c|} \hline K_1&K_2&K_3\\ \hline \left( \begin{array}{crcr} 1&1&1&1\\ 1&1&1&1\\
1&-1&1&-1\\1&-1&1&-1\end{array} \right)&
\left(\begin{array}{cccccccc}1&1&1&1&1&1&1&1\\1&1&1&1&1&1&1&1\\1&1&1&1&1&1&1&1\\1&1&1&1&1&1&1&1\\1&1&-&-&1&1&-&-\\1&1&-&-&1&1&-&-\\
1&1&-&-&1&1&-&-\\1&1&-&-&1&1&-&- \end{array}\right) &\left(\begin{array}{cccccccc}1&1&1&1&1&1&1&1\\1&1&1&1&1&1&1&1\\1&1&1&1&1&1&1&1\\1&1&1&1&1&1&1&1\\
1&-&1&-&1&-&1&-\\1&-&1&-&1&-&1&-\\
1&-&1&-&1&-&1&-\\1&-&1&-&1&-&1&- \end{array}\right)
\\ \hline
\end{array}$$

This basis may be extended to a full basis for 2-cocycles over $G_2^t$, by simply juxtaposing a basis for 2-coboundaries over $G_2^t$. Such a basis
was described in \cite{AAFR07} (according to
Lemma \ref{basecobordes}), $$\begin{tabular}{|c|l|} \hline $t$& basis for 2-coboundaries \\
\hline $\mbox{[}t]_2=1$&$\langle
\partial_2,\ldots,\partial_{4t-2}\rangle$ \\ \hline $\mbox{[}t]_2=0$ & $\langle
\partial_2,\ldots,\partial_{4t-3}\rangle$\\ \hline \end{tabular}$$

\begin{remark} All the generators above coincide with those of \cite{AAFR07} (notice that $G_2^t$ was denoted $G_4^t$ there), excepting
$BN_t \otimes 1_4$, which is substituted by $\beta _3=1_q\otimes BN_{2^r}\otimes 1_4$, for $t=2^rq$. It may be checked that $$\D
\beta_3=BN_t \otimes 1_4\cdot \prod _{k=0}^{\lfloor \frac{q}{2}\rfloor-1} \partial _{k2^{r+3}+2^{r+2}+1} \ldots \partial _{k2^{r+3}+2^{r+3}}.$$
\end{remark}

The cohomological reduction method also computes $H^3(cohG_2^t;\Z_2)=\left\{ \begin{array}{ll} \Z_2^4 & \mbox{ for $t$ odd}\\  \Z_2^{10} & \mbox{ for $t$ even} \end{array}\right.$ and a basis $\{\psi_1, \ldots , \psi _{10}\}$  (just $\{\psi_7,\ldots, \psi _{10}\}$ if $t$ is odd) for 3-cocycles over $G_2^t$  is given by:

 \begin{itemize}

 \item $\psi_1=\{ J_{4t},J_{4t},J_{4t},J_{4t},BN_t \otimes 1_4,BN_t \otimes 1_4,BN_t \otimes 1_4,BN_t \otimes 1_4, \stackrel{\frac{t}{2}}{\ldots}, J_{4t},J_{4t},J_{4t},J_{4t},$ $BN_t \otimes 1_4,BN_t \otimes 1_4,BN_t \otimes 1_4,BN_t \otimes 1_4 \}.$

 \item $\psi_2=\{ J_{4t},J_{4t},BN_t \otimes 1_4,BN_t \otimes 1_4, \stackrel{t}{\ldots}, J_{4t},J_{4t},BN_t \otimes 1_4,BN_t \otimes 1_4 \}.$

 \item $\psi_3=\{ J_{4t},BN_t \otimes 1_4, \stackrel{2t}{\ldots} ,J_{4t},BN_t \otimes 1_4\}.$

  \item $\psi_4=\{ J_{4t},J_{4t},1_{\frac{t}{2}}\otimes A \otimes 1_2,1_{\frac{t}{2}}\otimes A \otimes 1_2, \stackrel{t}{\ldots},J_{4t},J_{4t},1_{\frac{t}{2}}\otimes A \otimes 1_2,1_{\frac{t}{2}}\otimes A \otimes 1_2 \}$, for $A= \left( \begin{array}{cccc} 1&1&1&1\\1&1&1&1\\1&-&1&-\\1&-&1&- \end{array}\right)$.

\item $\psi_5=\{ J_{4t},1_{\frac{t}{2}}\otimes A \otimes 1_2, \stackrel{2t}{\ldots}, J_{4t},1_{\frac{t}{2}}\otimes A \otimes 1_2 \}$.

\item $\psi_6=\{ J_{4t},1_{\frac{t}{2}}\otimes \left( \begin{array}{cc}  J_4 &J_4\\ B&B \end{array}\right), \stackrel{2t}{\ldots} , J_{4t},1_{\frac{t}{2}}\otimes \left( \begin{array}{cc}  J_4 &J_4\\ B&B \end{array}\right)\}$, for $B= \left( \begin{array}{cccc} 1&-&1&-\\1&-&1&-\\1&-&1&-\\1&-&1&- \end{array}\right)$.

\item $\psi_7=\{ J_{4t},J_{4t},1_t \otimes BN_2 \otimes 1_2,1_t \otimes BN_2 \otimes 1_2, \stackrel{t}{\ldots},J_{4t},J_{4t},1_t \otimes BN_2 \otimes 1_2,1_t \otimes BN_2 \otimes 1_2 \}.$

\item $\psi_8=\{ J_{4t},1_t \otimes BN_2 \otimes 1_2, \stackrel{2t}{\ldots} , J_{4t},1_t \otimes BN_2 \otimes 1_2\}.$

\item $\psi_9=\{ J_{4t},1_t \otimes A, \stackrel{2t}{\ldots}, J_{4t},1_t \otimes A \}.$

\item $\psi_{10}=\{ J_{4t},1_{2t} \otimes BN_2, \stackrel{2t}{\ldots} ,J_{4t},1_{2t} \otimes BN_2\}.$

\end{itemize}

This basis may be extended to a full basis for 3-cocycles over $G_2^t$, by simply juxtaposing a basis for 3-coboundaries over $G_3^t$. According to
Lemma \ref{basecobordes}, such a basis is given by $$\begin{tabular}{|c|l|} \hline $t$& basis for 3-coboundaries \\
\hline $\mbox{[}t]_2=1$&$\langle
\partial_1,\ldots,\partial_{16t^2-4t-2},\partial_{16t^2-4t+1}\rangle$ \\ \hline $\mbox{[}t]_2=0$ & $\langle
\partial_1,\ldots,\partial_{16t^2-8t-3},\partial_{16t^2-8t+1},\ldots,\partial_{16t^2-4t-3},\partial_{16t^2-4t+1},\ldots,\partial_{16t^2-4t+3}\rangle$\\
\hline \end{tabular}$$
 by direct inspection.

\subsection{Basis for 3-cocycles over $G_3^t=\Z_{2t}$}$\mbox{ }$\\

Consider the family $G_3^t=\Z _{2t}$. We now apply the cohomological reduction method for calculating a basis for $n$-cocycles over $G$. It may be checked (see \cite{EM54}) that a homological model for $G_3^t$ is given by ${\cal B}_2^*=\{e_1^*\}$, ${\cal B}_3^*=\{v_1^* \}$, ${\cal B}_4^*=\{w_1^* \}$,
$d^2=0$, $d^3=0$, $F[i|j|k]=[ k[i+j\geq 2t]]_2\cdot v_1$.

From these data, it may be checked that $\begin{array}{|c|c|c|c|} \hline D_2&Q^{-1}_3&D_3&P_3\\ \hline  0 & 1 & 0 & 1 \\
\hline
\end{array}$.

Thus $H^3(cohG_3^t;\Z_2)=\Z_2$ and ${\cal B}_3^*=\{v_1^* \}$ is a basis for representative 3-cocycles over $cohG_3^t$.
Accordingly, a basis for representative 3-cocycles over $G_3^t$ is given by $\{F^*(v_1^*)\}$, whose horizontal sections are given by
%
%
$\{J_{2t},BN_{2t}, \stackrel{t}{\ldots},J_{4t},BN_{2t} \}$.
%

This basis may be extended to a full basis for 3-cocycles over $G_3^t$, by simply juxtaposing a basis for 3-coboundaries over $G_3^t$, which is given by $\langle
\partial_1,\ldots,\partial_{4t^2-2t}\rangle$ by direct inspection.

\subsection{Basis for 3-cocycles over $G_4^t=D_{2t}$, $t$ odd}$\mbox{ }$\\

Unfortunately, it is known that $H^3 (D_{2t})=0$ for odd $t$, and therefore a basis for 2-cocycles over $G_4^t$ consists in a basis for 2-coboundaries over $G_4^t$, such as $\langle \partial_2,\ldots,\partial_{2t+1},\partial_{4t},\ldots ,\partial_{2t^2+1},\partial_{2t^2+3},\ldots , \partial_{2t^2+t}, \partial_{2t^+t+2},\ldots, \partial_{4t^2}\rangle$.

\subsection{Calculating 3-dimensional Hadamard matrices}

Using the basis for 3-cocycles over $G_1^1=G_2^1=\Z_2^2$ calculated before, we have performed an exhaustive search for $3$-dimensional 3-cocyclic Hadamard matrices over $\Z_2^2$. This search yields that there are 64 improper such matrices, none of which is in addition proper.  For instance, the matrix related to the product of 3-coboundaries $\partial _4 \partial _7 \partial _{10} \partial _{13}$, whose horizontal sections are given by $$\{ \left( \begin{array}{cccc}1&1&1&1\\ 1&1&-&-\\ 1&-&1&-\\ -&1&1&-\end{array} \right), \left( \begin{array}{cccc} 1&1&-&-\\1&1&1&1\\ 1&-&-&1\\ -&1&-&1 \end{array} \right), \left( \begin{array}{cccc} 1&-&1&-\\ 1&-&-&1\\ 1&1&1&1\\ -&-&1&1 \end{array} \right), \left( \begin{array}{cccc} -&1&1&- \\ -&1&-&1\\ -&-&1&1\\ 1&1&1&1 \end{array} \right) \}.$$

An exhaustive computer search for $3$-dimensional Hadamard matrices 3-cocyclic over $G_3^2=\Z_4$, yields that there are 32 improper such matrices, none of which is in addition proper. For instance, the matrix related to the product of 3-coboundaries $\partial _4 \partial _7 \partial _8 \partial _9$, whose horizontal sections are given by $$\{ \left( \begin{array}{cccc}1&1&1&1\\ 1&-&-&1\\ -&1&-&1\\ 1&1&-&-\end{array} \right), \left( \begin{array}{cccc} 1&1&-&-\\1&-&1&-\\ -&1&1&-\\ 1&1&1&1 \end{array} \right), \left( \begin{array}{cccc} 1&-&1&-\\ 1&1&-&-\\ -&-&-&-\\ 1&-&-&1 \end{array} \right), \left( \begin{array}{cccc} -&1&1&- \\ -&-&-&-\\ 1&1&-&-\\ -&1&-&1 \end{array} \right) \}.$$

Since there is more interest in improper $n$-dimensional Hadamard matrices of even order $v=4k+2$ (in particular, for $n=3$ and $v \neq 2 \cdot 3^b$),  we have also performed a partial heuristic search for 3-dimensional improper Hadamard matrices of order $2t$ over $G_3^t$ and $G_4^t$, for $t=3,5$, from which unfortunately we have not got any 3-dimensional 3-cocyclic improper Hadamard matrix. Contributions in this sense (with these or other groups of order $v=4k+2$) would be appreciated, since improper $3$-dimensional Hadamard matrices of order different to $2 \cdot 3^b$ are still to be discovered.

\section{Conclusions and further work}

It is well-known that there exists a proper $n$-dimensional Hadamard matrix of order $v\equiv 0\; \mbox{mod 4}$ for every $n\geq 3$ if and only if there exists a planar Hadamard matrix of order $v$. However,  few is known about the existence of improper $n$-dimensional Hadamard matrices of order $v$. It is evident that $v$ must be even, but surprisingly it need not to be a multiple of 4. And it is not known whether improper $n$-dimensional Hadamard matrices exist for all even orders $v=2t$, $n \geq 3$, for odd $t$. Furthermore, only $3$-dimensional improper Hadamard matrices of order $2 \cdot 3^b$ are known, $b\geq 1$ (see \cite{Yan01} for more details).

In this paper we have been concerned with higher dimensional Hadamard matrices, and the way in which the cocyclic framework may be introduced in this context, going beyond de Launey's works in \cite{dLa89,dLH93,dLF11}.

In particular,  we have  provided a method for constructing $n$-cocyclic $n$-dimensional matrices over some finite groups $G$, from which a deeper search for improper/proper Hadamard matrices might be performed. The input data of our process is a group $G$ for which a cohomological model is known, so that the cohomological reduction method may be straightforwardly applied.

Focusing in the case $n=2$ (which is of special interest, e.g. for looking for cocyclic Hadamard matrices), it follows that there is no need to distinguish between inflation and transgression 2-cocycles, as it has traditionally been the case. Some examples have been given in Section 3.

We have also provided some examples of $3$-dimensional 3-cocyclic Hadamard matrices, progressing on the basis for $3$-cocycles provided by our algorithm.

Finally, we would like to conclude this paper proposing several problems concerning higher dimensional (cocyclic) Hadamard matrices, which should be studied in a near future, such as:

\begin{enumerate}

\item Characterize {\em $n$-cocyclic} $n$-dimensional proper and improper Hadamard matrices.

One could think that an $n$-dimensional $n$-cocyclic matrix would consist of $2$-cocyclic layers. This is not true at all, as it will be discussed
elsewhere. Thus the traditional cocyclic test for 2-dimensional 2-cocyclic Hadamard matrices cannot be naturally extended to the
$n$-dimensional case so far. A deeper analysis must be done.

\item Look for $n$-dimensional proper and improper $n$-cocyclic Hadamard matrices.

This is a very difficult task, since the search space seems to grow drastically in exponential size (e.g., the basis for 3-dimensional 3-cocyclic
matrices over $G_1^t$ and $G_2^t$, $t$ even, consists of $16t^2-4t+3$ and $16t^2-4t+7$ generators, respectively). 

Maybe one should think of constructing $n$-dimensional improper Hadamard matrices based on planar cocyclic matrices satisfying some certain constraints, in light of de Launey's fruitful way for constructing $n$-dimensional proper Hadamard matrices from $2$-cocyclic Hadamard matrices, as described in \cite{dLa89,dLH93,dLF11}.



\item Determine whether any of the already known construction methods for generating higher dimensional proper Hadamard matrices from
2-dimensional ones (see \cite{Hor07,Yan01} for details) involves $n$-cocyclic matrices.

Furthermore, is it possible to derive  a method for constructing $n$-dimensional proper/improper $n$-cocyclic Hadamard matrices from pro\-per/improper $k$-cocyclic Hadamard matrices, $k<n$ (e.g. via the cup product in cohomology)?

\item Do a $n$-dimensional improper ($n$-cocyclic) Hadamard conjecture make sense?

Proper higher dimensional Hadamard matrices might exist only for orders 1,2 and a multiple of 4, and they do exist if and only if a planar Hadamard matrix of the same order do exist. This reduces the problem to the Hadamard Conjecture. Nevertheless, as it has been discussed earlier, only higher dimensional improper Hadamard matrices of order $2 \cdot 3^b$ are known. So the question is widely open for the improper case.

\end{enumerate}

\vspace*{.5cm}
\noindent{\bf \large Acknowledgements}

\vspace*{.5cm}

The authors want to express their gratitude to the anonymous referees for their valuable advices and suggestions, which have helped to
improve the readability of the paper for a better understanding.

\end{document}